\documentclass[12pt,a4paper]{article}
\usepackage{amsmath}
\usepackage{amssymb}
\usepackage[english]{babel}
\usepackage[left=1cm,right=1cm,top=2cm,bottom=2cm]{geometry}
\usepackage{color}
\usepackage{tikz}
\usepackage{abstract}
\usetikzlibrary{calc,intersections}
\newenvironment{nscenter}
 {\parskip=0pt\par\nopagebreak\centering}
 {\par\noindent\ignorespacesafterend}

\begin{document}

\begin{center}\textbf{Golubeva V. A., Ivanov A. N.}\end{center}

\begin{center}
{\textbf{Fuchsian systems for Dotsenko--Fateev multipoint correlation functions and similar integrals of hypergeometric type}}
\end{center}

\begin{abstract}
\fontdimen2\font=10pt
\noindent The Dotsenko--Fateev integral {is an analytic function of one complex variable expressing the amplitude in the 4-point correlator of the 2D conformal field theory. Dotsenko--Fateev found ODE of third order with Fuchsian singularities satisfied by their integral. In the present paper, this work is extended to generalized Dotsenko--Fateev integrals, in particular those associated to arbitrary multipoint correlators, and Pfaff systems of PDE of Fuchsian type are constructed for them. The ubiquity of the Fuchsian systems is in that they permit to obtain
local} expansions  of solutions in the {neighbor\-hoods} of singularities of the system.

\vspace{1mm}\noindent\textit{Keywords}: Dotsenko--Fateev integral, systems of partial differential equations, Pfaffian systems of Fuchsian type, hyperplane arrangements, hypergeometric
functions
\end{abstract}

\begin{nscenter}\textbf{\hspace{0.5mm}1. Introduction}\end{nscenter}

\vspace{1mm} In two-dimensional conformal field theory{,} a $(n+2)$-point correlation { function  containing} a { third order conformal} operator {  can be expressed, using the projective invariance to fix arbitrary three points, in terms of} special correlators of the following form{:}

$${ \langle}\phi_{\alpha_{k_0,l_0}}(0)\phi_{\alpha_{k_1,l_1}}(1)\phi_{\alpha_{1,3}}(x_2)...\phi_{\alpha_{1,3}}(x_{n})\phi_{\alpha_{k_3,l_3}}(
\infty){ \rangle=}$$$$\int\limits_{C_1 \times C_2} { \langle}V_{\alpha_{k_1,l_1}}(0)V_{\alpha_{k_2,l_2}}(1)V_{\alpha_{1,3}}(x_2)...V_{\alpha_{1,3}}(x_{n})V_{\alpha_{k_3,l_3}}(
\infty)V_{\alpha_{+}}(t)V_{\alpha_{+}}(s){ \rangle \, dt ds  \  \sim} \eqno(1.1)$$$$ \prod_{i=2}^{n} x_i^{2\alpha_{k_1,l_1}\alpha_{1,3}}\prod_{i=2}^{n} (1-x_i)^{2\alpha_{k_2,l_2} \alpha_{1,3}} \int\limits_{C_1
\times C_2} t^{2\alpha_{k_1,l_1}\alpha_{+}} s^{2\alpha_{k_1,l_1}\alpha_{+}} (t-1)^{2\alpha_{k_2,l_2}\alpha_{+}} (s-1)^{2\alpha_{k_2,l_2}\alpha_{+}}{ \times}$$$$\prod_{i=2}^{n}
(t-x_i)^{2\alpha_{1,3}\alpha_{+}}(s-x_i)^{2\alpha_{1,3} \alpha_{+}} (t-s)^{2\alpha_{+}^2} dt ds,$$

\noindent where $\alpha_{k,l}=\frac{1-k}{n} \alpha_{-}+\frac{1-l}{n} \alpha_{+}$ {are quantum} charges, $\alpha_{\pm}$ are two roots of the equation $\Delta_{\alpha_{\pm}}=\alpha^2_{\pm}-2\alpha_{\pm}\alpha_0=0$ ({i.e.} the
conformal dimension of the integral operator $\int\limits_{C_i} V_{\alpha_{\pm}}(t) dt$ should {be zero}; $\alpha_0$ is some constant), $V_{\alpha}$ is the {Coulomb operator with charge} $\alpha$ ( {i.e.} the exponent
$e^{i\alpha\phi}$ of the free field $\phi$), $C_i$ are the integration paths satisfying the conditions $\overset{\circ}{C_i} \subset \bar{\mathbb{C}} \setminus \{0, 1, x_2,...,x_n, \infty\},~
\partial C_i \in \{0, 1, x_2, ..., x_n, \infty\}$. The integral {contained} in the above expression for the case $n=2$ {is} called Dotsenko{--Fateev} integral. {It} was investigated in [1, 2]{,  where an} ordinary differential equation of third order was found{, of which it is a solution.}

The present paper is devoted to the investigation of the integral defining the correlation function (1.1) for arbitrary $n$, that is the integral of the following form{:}

$$I_n(a_0,a_1,..., a_n, g; x_2, ..., x_n)=\int\limits_{C_1 \times C_2} t^{a_0} s^{a_0} (t-1)^{a_1} (s-1)^{a_1} \prod\limits_{i=2}^n (t-x_i)^{a_i}(s-x_i)^{a_i} (t-s)^g  dt ds, \eqno(1.2)$$

\noindent where $a_i \in \mathbb{C}$  {are} complex parameters, {and} $C_j$ {are} integration paths satisfying the conditions $\overset{\circ}{C_j} \subset \bar{\mathbb{C}} \setminus \{0, 1,
x_2,...,x_n, \infty\},~ \partial C_j \in \{0, 1, x_2, ..., x_n, \infty\}$. {Moreover, the} natural generalization of the Dotsenko{--}Fateev integral

$$J_n(a, b, c, g; x_1,...,x_n)=\int\limits_{C_1 \times ... \times C_n} \prod\limits_{i=1}^{n}t_{i}^{a} \prod\limits_{i=1}^{n}(t_i-1)^{b} \prod\limits_{i=1}^{n}(t_i-x_i)^{c}
\prod\limits_{i<j}(t_i-t_j)^{g} dt_1...{dt_n} \eqno(1.3)$$

\noindent {is studied,} where the integration paths $C_i$ satisfy the same conditions. {Fuchsian} systems ({i.e. Pfaff systems} of Fuchsian type) { are obtained, which admit certain vector-valued functions associated to $I_n$, $J_2$ as solutions. For $J_2$, two Fuchsian systems are given}, one of them {being} equivalent to the Dotsenko{--}Fateev {third order} equation. The derivation of the Fuchsian systems is based on the theory of hyperplane arrangements{, which provides some tools for handling integrals of hypergeometric type, whose integrands are products of complex powers of linear polynomials}. The Fuchsian systems are important for applications{,} because they {permit} to obtain local expansions  of solutions in the neighborhoods of singularities of the system [10]. Besides, {a} system of partial differential equations {for $J_n$ is} obtained. This system is {analogous} to the Appell and {Kamp\'e de F\'eriet} equations for hypergeometric functions of two complex variables $F_1, F_2, F_3, F_4$.

The structure of the present paper {is} the following. In {Sec. 2, a} detailed description of the method {for deriving Fuchsian} systems for the vector-valued functions associated {to} hypergeometric integrals is given. This method {consists} in {a} description of {a} basis of sections {of the} cohomology bundle corresponding to the {integral under consideration, followed by a} computation of certain connection matrix in this basis. In {Sec. 3, an} explicit description of {the} basis is given. In {Sec. 4, some auxiliary} formulas{, necessary for the} computation of {the} connection matrix are {obtained}. In {Sec. 5, 6,} the cases of the integrals $I_n, J_2$ are considered. In {Sec. 7, a} system of partial differential {equations for $J_n$ is obtained}. At last, in {Sec. 8, another Fuchsian system for the Dotsenko--Fateev integral, equivalent to the Dotsenko--Fateev third order equation, is obtained by elementary methods}.\\

\vspace{2mm}\begin{nscenter}\textbf{\hspace{0.5mm}2. {Method for deriving Fuchsian} systems}\end{nscenter}

\vspace{1mm} Let us consider the following hypergeometric integral{:}

$$I(x_1,...,x_n)=\int\limits_{\sigma} \Phi dt, ~~ \Phi = \prod\limits_{i=1}^{l} \alpha_i (x_1,...,x_n;t_1,...,t_m)^{\lambda_i} , ~~ dt=dt_1 ... dt_m, \eqno(2.1)$$

\noindent where $\alpha_i (x_1,...,x_n;t_1,...,t_m)=a_{i0}(x_1,...,x_n)+\sum\limits_{j=1}^{l} a_{ij}(x_1,...,x_n) t_j$, $a_{ij}(x_1,...,x_n)$ are polynomials of variables $x_1,...,x_n$; $\lambda_i
\in \mathbb{C}$ are complex parameters with sufficiently large real parts, $\sigma$ is the integration domain such that $\Phi|_{\partial\sigma} \equiv 0$, $\sigma$ ($\sigma$ may depend on the variables $x_1, ..., x_n$).

Let $Z=\big\{(x,t) \in \mathbb{C}^n \times \mathbb{C}^m ~ | ~ \Phi(x,t)=0 \big\}${, and let} $\text{pr}: Z \longrightarrow \mathbb{C}^n$ be the corresponding {projection.} Then we have the family of hyperplane arrangements $Z_x=\text{pr}^{-1}(\{x\}),~ x \in \mathbb{C}^n$. We denote the {complement} of the arrangement $Z_x$ in $\mathbb{C}^m$ by $\mathbb{C}_{x}^{m}=\mathbb{C}^m \setminus Z_x$. For each $x \in \mathbb{C}^n${,} we consider the operator $\nabla_x: \Omega^{\boldsymbol{\cdot}}(*Z_x) \longrightarrow \Omega^{\boldsymbol{\cdot}+1}(*Z_x)$ defined by the following condition{:}

$$d(\Phi \psi)=\Phi \nabla_x \psi ~~~ \forall x \in \mathbb{C}^n ~~~ \forall \psi \in \Omega^{\boldsymbol{\cdot}}(*Z_x),$$

\noindent where $\Omega^{k}(*Z_x)$ is the vector space of rational {differential} $k$-forms on $\mathbb{C}^m$ with poles along $Z_x$, $d$ is the differential with respect to the variables $t_1,...,t_m$. It {is} easy to see that this operator has the form $\nabla_x=d+\omega_{\lambda} \wedge$, where $\omega_{\lambda}=\sum\limits_{i=1}^{l} \lambda_i \omega_i \wedge,~ \omega_i=d \ln\alpha_i \in \Omega^{1}(*Z_x)$. Besides, this operator is the differential for the following complex

$$0 \longrightarrow \Omega^{0}(*Z_x) \overset{\nabla_x}{\longrightarrow} \Omega^{1}(*Z_x) \overset{\nabla_x}{\longrightarrow} ... \overset{\nabla_x}{\longrightarrow} \Omega^{m}(*Z_x)
\longrightarrow 0. \eqno(2.2)$$

\noindent The value of the integral $\int\limits_{\sigma} \Phi \psi$ at the point $x \in \mathbb{C}^n$ depends on the class of the differential form $\psi$ in the highest cohomology $\text{H}^m(\Omega^{\boldsymbol{\cdot}}(*Z_x),\nabla_x)$ by virtue of the following {identities:}

$$\int\limits_{\sigma} \Phi \nabla_x \psi=\int\limits_{\sigma} d(\Phi \psi)=\int\limits_{\partial\sigma} \Phi \psi=0 \eqno(2.3)$$

\noindent whenever $\Phi$ decreases rapidly enough near $\partial\sigma$, so that the integrals converge.

Further consider the operator $\nabla': \Omega^{0,\boldsymbol{\cdot}}(*Z) \longrightarrow \Omega^{1,\boldsymbol{\cdot}}(*Z)$ defined by the condition

$$d'(\Phi \psi)=\Phi \nabla' \psi ~~~ \forall \psi \in \Omega^{0,\boldsymbol{\cdot}}(*Z),$$

\noindent where $\Omega^{p,q}(*Z)$ is the vector space of rational differential $(p+q)$-forms on $\mathbb{C}^n \times \mathbb{C}^m$ {with poles along $Z$,} which are $q$-forms with respect to the variables $t_1,...,t_m$; $d'$ is the {differential} with respect to the variables $x_1, ..., x_n$. Obviously{, the} operator $\nabla'$ has the expression $\nabla'=d'+\sum\limits_{i=1}^{l} \lambda_i \omega'_i \wedge$, where $\omega'_i=d'
\ln\alpha_i \in \Omega^{1,0}(*Z)$.

Now {fix some} point $x_0 \in \mathbb{C}^n$ and consider {the} subset $X\subset \mathbb{C}^n$ { of all $x \in X$ such that the arrangements $Z_x$ are combinatorially equivalent to $Z_{x_0}$} (we restrict {ourselves to an open set $X=\mathbb{C}^n \setminus \cup L_i$, where $L_i$ are algebraic subvarieties of $\mathbb{C}^n$}). Then $\mathcal{H}=\underset{x \in
X}{\bigcup} \text{H}^m(\Omega^{\boldsymbol{\cdot}}(*Z_x),\nabla_x)$ can be naturally endowed {with a} structure of {a} vector bundle over $X${,} and the operator $\nabla'$ can be considered as {a} connection {on} this vector bundle. Suppose that the collection of logarithmic differential forms $\{ \eta_i \in  \Omega^{m,0}(*Z)\}$ defines the basis of global sections ({i.e.} it is {a} fiberwise basis) of the vector bundle $\mathcal{H}$. Then the action of the operator $\nabla'$ {on the} forms $\eta_i$ is {characterized} by the formula

$$\nabla'\eta_i=\sum\limits_{i=1}^{n} dx_i \wedge \left(\sum\limits_{j} f_{ij}(x_1,..,x_n) \eta_j+\nabla(\psi_i)\right). \eqno(2.4)$$

\noindent {Upon integration over $\sigma$,} we obtain

$$d'\left(\int\limits_{\sigma} \Phi \eta_i\right)=\int\limits_{\sigma} d'(\Phi \eta_i)=\int\limits_{\sigma} \Phi \nabla'(\eta_i){=}$$$$\sum\limits_{i=1}^{n} dx_i \wedge \left(\sum\limits_{j}
f_{ij}(x_1,..,x_n) \int\limits_{\sigma} \Phi \eta_j\right)=\sum\limits_{j} \Omega_{ij} \wedge \int\limits_{\sigma} \Phi \eta_j, \eqno(2.5)$$

\noindent where $\Omega_{ij}$  {are differential} 1-forms on $X$ (in the case if $\sigma$ doesn't depend on $x$, the first equality of (2.5) is evident; in the general case this equality follows from the condition $\Phi|_{\partial \sigma} \equiv 0$). Denote $\left(\int\limits_{\sigma} \Phi \eta_i\right)^T$ by $f$, then the relations (2.5) {have} the form of the Pfaffian system $d'f=\Omega f$ on $X$. Moreover, if we choose all bounded chambers in $\mathbb{R}^m_{x} \subset \mathbb{C}^m_{x}$ as integration domains $\sigma_j$ we obtain the basis of solutions $\left(\int\limits_{\sigma_j} \Phi \eta_i\right)^T$ of this Fuchsian system [3].

Therefore, the problem of construction of Pfaffian system {is} reduced to the {construction} of {a} collection of rational logarithmic differential forms on $\mathbb{C}^n \times \mathbb{C}^m$ {forming a} basis of global sections of the vector bundle $\mathcal{H}$, and the computation of {the} connection matrix $\Omega$ in this basis. These problems are solved {by methods of the} theory of hyperplane arrangements which will be described {in} the next section. {We will provide a} collection of {$\eta_i$} such that the corresponding Pfaffian system is Fuchsian, {the} connection matrix $\Omega$ {having} the form $\sum\limits_{i} A_i \frac{d' L_i}{L_i}$ with constant matrices $A_i$.\\

\begin{nscenter}\textbf{\hspace{0.5mm}3. {A} basis of global sections }\end{nscenter}

\vspace{1mm} At first we describe {a} basis of the vector space $\text{H}^m(\Omega^{\boldsymbol{\cdot}}(*\mathcal{A}),\nabla)$, where $\mathcal{A}=Z_{x_0}, \nabla=\nabla_{x_0}$, and then extend it to {the} desired collection of differential forms on $X$. Consider the graded subalgebra $B^{\boldsymbol{\cdot}}(\mathcal{A})$ of algebra $\Omega^{\boldsymbol{\cdot}}(*\mathcal{A})$ generated by {the} differential forms $\omega_H=d \ln \alpha_H$ for all $H \in \mathcal{A}$. The inclusion $B^{\boldsymbol{\cdot}}(\mathcal{A}) \subset \Omega^{\boldsymbol{\cdot}}(*\mathcal{A})$ {induces an} isomorphism of corresponding cohomologies $\text{H}^{\boldsymbol{\cdot}}(B^{\boldsymbol{\cdot}}(\mathcal{A}),\omega_\lambda)
\simeq \text{H}^{\boldsymbol{\cdot}}(\Omega^{\boldsymbol{\cdot}}(*\mathcal{A}),\nabla)$ (for a proof, see [7]).

Moreover the algebra $B^{\boldsymbol{\cdot}}(\mathcal{A})$ is isomorphic to the graded Orlik{--}Solomon algebra $A(\mathcal{A})$ of the arrangement $\mathcal{A}$. This algebra $A(\mathcal{A})$ is defined as the quotient of the {Grassmannian} algebra $\Lambda\left(\oplus_{H\in\mathcal{A}}\mathbb{C}e_H\right)$, where $e_H$ are the formal variables corresponding to all hyperplanes of arrangement $\mathcal{A}$, by the ideal generated be the elements of the form $e_{H_{i_1}} ... e_{H_{i_m}}$ for $H_{i_1} \cap ... \cap H_{i_m} =\emptyset$ and the elements of the form $\partial(e_{H_{i_1}} ... e_{H_{i_m}})=\sum\limits_{k=1}^{m}(-1)^{k-1} e_{H_1}...e_{H_{k-1}} e_{H_{k+1}}...e_{H_m}$ for $H_{i_1} \cap ... \cap H_{i_m} \neq \emptyset, ~ \text{codim}(H_{i_1} \cap ... \cap H_{i_m})< m$. This isomorphism $B^{\boldsymbol{\cdot}}(\mathcal{A}) \simeq A(\mathcal{A})$ is induced by the map $a_H \longmapsto \omega_H$ where $a_H$ is the projection of $e_H$ in $A(\mathcal{A})$ (for a proof, see [9]).

Fix further a linear order on the set of hyperplanes from the arrangement $\mathcal{A}$ and consider the set $2^{\mathcal{A}}$ of all subsets $S=\{H_{i_1},...,H_{i_k}\}$ of the arrangements $\mathcal{A}$ ({final} results {do} not depend on the choice of the order). We call {a} subset $S$ dependent, if $\cap S=H_{i_1} \cap ... \cap H_{i_m} {\neq} \emptyset$, $\text{codim}(\cap S) < |S|$, and independent, if $\cap S \neq \emptyset$, $\text{codim}(\cap S) = |S|$. {A dependent subset of $2^{\mathcal{A}}$, minimal by inclusion,} is called {a} circuit. Also we call { $S$} a broken circuit, if there exists {a} hyperplane $H \in \mathcal{A}$, such that $H < \text{min}(S)$ and $S \cup {H}$ is a circuit. Now consider the collection $\text{nbc}(\mathcal{A})$ of subsets $S \in
2^{\mathcal{A}}$, such that $\cap S \neq \emptyset$ and $S$ {contains no} broken circuits. It {is} easy to see that $\text{nbc}(\mathcal{A})$ is closed with respect to taking subsets of { $S$}. Hence, $\text{nbc}(\mathcal{A})$ is {an} abstract simplicial complex (the dimension of {a} simplex $S$ is equal {to} $|S|-1$). Further we will suppose that the dimension of $\text{nbc}(\mathcal{A})$ is equal {to} $m-1$.

Now we consider cochain complex $C^{\boldsymbol{\cdot}}(\text{nbc}(\mathcal{A}),\mathbb{C})$ of the simplicial complex $\text{nbc}(\mathcal{A})$ with coefficients in $\mathbb{C}$. The map $\Theta:
C^{\boldsymbol{\cdot}-1}(\text{nbc}(\mathcal{A}),\mathbb{C}) \longrightarrow A^{\boldsymbol{\cdot}}(\mathcal{A})$ defined by the {formula}

$$\Theta(\alpha)= {{\scalebox
                {1.4}
                {${\displaystyle\sum\limits_{S \in \text{nbc}(\mathcal{A}), |S|=q}}$}}} \alpha(S) {
                {\scalebox
                {1.2}
                {${\displaystyle\bigwedge_{p=1}^{q}}$}}}\left(\sum\limits_{\cap_{k=p}^{q} H_{i_k} \subset H} \lambda_H a_H \right) \eqno(3.1)$$

\noindent  {induces an} isomorphism of cochain complexes $C^{\boldsymbol{\cdot}-1}(\text{nbc}(\mathcal{A}),\mathbb{C}) \simeq (A^{\boldsymbol{\cdot}}(\mathcal{A}),a_{\lambda}\wedge)${,} where $a_{\lambda} \in A^{\boldsymbol{\cdot}}(\mathcal{A})$ {correspond} to the differential {forms} $\omega_{\lambda} \in B^{\boldsymbol{\cdot}}(\mathcal{A})$ (for a proof, see [3]). By this reason we have the following sequence {of} isomorphisms

$$\text{H}^{q}(\Omega(*\mathcal{A}),\nabla) \simeq \text{H}^{q}(B^{\boldsymbol{\cdot}}(\mathcal{A}),\omega_\lambda) \simeq \text{H}^{q}(A^{\boldsymbol{\cdot}}(\mathcal{A}), a_{\lambda}\wedge) \simeq
\text{H}^{q-1}(\text{nbc}(\mathcal{A}),\mathbb{C}). \eqno(3.2)$$

We can {obtain a} basis of the highest cohomology $\text{H}^{m-1}(\text{nbc}(\mathcal{A}),\mathbb{C})$ of the simplicial complex $\text{nbc}(\mathcal{A})$ {in} the following way. Choose from the highest dimension simplices of {$\text{nbc}(\mathcal{A})$ a subset  $\beta\text{nbc}(\mathcal{A})$ of} simplices $S$ satisfying the {following} condition: for all {hyperplanes} $H \in S${,} there exists {a} hyperplane $H' \in \mathcal{A}$,
such that $(S \setminus {H}) \cup {H'}$ {is a maximal independent} subset of $2^{\mathcal{A}}$. Then the dual elements of $\text{H}^{m-1}(\text{nbc}(\mathcal{A}),\mathbb{C})$ corresponding to $\beta \text{nbc}(\mathcal{A})$ form {a} basis (for a proof, see [8]).

Now{, using} (3.2){, we obtain a} basis of $\text{H}^{m}(\Omega(*\mathcal{A}),\nabla)$ {defined as follows:}

$$\theta(S)={
                {\scalebox
                {1.4}
                {${\displaystyle\bigwedge_{p=1}^{r}}$}}}\left(\sum\limits_{\cap_{k=p}^{q} H_{i_k} \subset H} \lambda_H d \ln \alpha_H \right), ~~ S \in \beta \text{nbc}(\mathcal{A}). \eqno(3.3)$$

\noindent {Applying this construction to our relative situation, where the hyperplanes of $\mathcal A=Z_x$ depend on $x$, we obtain a} collection of differential forms on $X$ {whose cohomology classes form a} basis of global sections of the vector bundle $\mathcal{H}$. {Such a} basis of global sections in the next sections will be called {a} $\beta \text{nbc}$-basis. \\

\begin{nscenter}\textbf{\hspace{0.5mm}4. The connection matrix in {a} $\beta\text{nbc}$-basis}\end{nscenter}

Let $T$ be the $(m+1)\times(l+1)$ matrix{, whose} first $l=|\mathcal{A}|$ columns are {filled with} the coefficients of linear polynomials $\alpha_{H_i},
i=1,...,|\mathcal{A}|=l$ (the last entries of these columns {being constant terms of the linear polynomials}) and the last column {is} $(0,0,...,1)^T$. For {an} ordered collection of indices $S=(i_1,..., i_{r})${,} we denote the collection $(i_1,...,i_{k-1},i_{k+1},...,i_{r})$ by $S_k$, the collection $(i_1,...,i_{r},j)$ by $(S,j)$ and we use the following notation{:}

$$\alpha_S= \alpha_{i_1}...\alpha_{i_r}, ~~~ \omega_S=\omega_{i_1} \wedge ... \wedge \omega_{i_r}, ~~~ \omega'_S=\sum\limits_{k=1}^{r}(-1)^{k+1} \omega'_{i_k} \wedge \omega_{S_k}, ~~~
\Delta_S=\det T_S, \eqno(4.2)$$

\noindent where $T_S$ is the {sub}matrix of {$T$ with columns} indexed by $S$. If $r=m+1${,} it is easy to see that $\omega_{S_k}=\frac{\alpha_{i_k}}{\alpha_{S}}\Delta_{(S_k,n+1)}dt$ and the following equalities are {true:}

$$\omega'_S=\sum\limits_{k=1}^{m+1}(-1)^{k+1} \omega'_{i_k} \wedge \omega_{S_k}=\frac{1}{\alpha_S}\sum\limits_{k=1}^{m+1}(-1)^{k+1} \Delta_{(T_k,n+1)} d_x \alpha_{i_k} \wedge
dt=\frac{1}{\alpha_S}\sum\limits_{k=1}^{m+1}(-1)^{k+1} d_x (\alpha_{i_k} \Delta_{(T_k,n+1)}) \wedge dt {-}$$$$\frac{1}{\alpha_S}\sum\limits_{k=1}^{l+1}(-1)^{k+1} \alpha_{i_k} d_x
\Delta_{(T_k,n+1)} \wedge dt=\frac{1}{\alpha_S} d_x\left(\sum\limits_{k=1}^{m+1}(-1)^{k+1} \alpha_{i_k} \Delta_{(T_k,n+1)}\right) \wedge dt {-}$$$$\sum\limits_{k=1}^{m+1}(-1)^{k+1} d_x \ln
\Delta_{(T_k,n+1)} \wedge \frac{\alpha_{i_k}}{\alpha_S} \Delta_{(T_k,n+1)} dt=\frac{1}{\alpha_S} d_x \Delta_{S} \wedge dt- \sum\limits_{k=1}^{m+1}(-1)^{k+1} d_x \ln \Delta_{(T_k,n+1)} \wedge
\omega_{S_k}{=}$$$$d_x \ln \Delta_{S} \wedge \frac{\Delta_S}{\alpha_S} dt- \sum\limits_{k=1}^{m+1}(-1)^{k+1} d_x \ln \Delta_{(T_k,n+1)} \wedge \omega_{S_k}{=}$$$$ d_x \ln \Delta_{S} \wedge
\left(\sum\limits_{k=1}^{m+1}(-1)^{k+1}\frac{\alpha_{i_k}}{\alpha_S} \Delta_{S_k} dt\right)- \sum\limits_{k=1}^{m+1}(-1)^{k+1} d_x \ln \Delta_{(T_k,n+1)} \wedge
\omega_{S_k}{=}$$$$\sum\limits_{k=1}^{m+1}(-1)^{k+1} d_x \ln \Delta_{S} \wedge \omega_{S_k} - \sum\limits_{k=1}^{m+1}(-1)^{k+1} d_x \ln \Delta_{(S_k,n+1)} \wedge
\omega_{S_k}=\sum\limits_{k=1}^{m+1}(-1)^{k+1} d_x \ln \frac{\Delta_S}{\Delta_{(S_k,n+1)}} \wedge \omega_{S_k}.$$

\noindent On the other hand, for arbitrary $r$ the following equalities {also hold:}
$$\nabla' \omega_S + \nabla \omega'_S=\sum\limits_{k=1}^{r}(-1)^{k+1} d_x \omega_{i_k} \wedge \omega_{S_k} + \sum\limits_{i=1}^{n}\lambda_i \omega'_{i} \wedge \omega_S +
\sum\limits_{k=1}^{r}(-1)^{k+1} d_t \omega'_{i_k} \wedge \omega_{S_k} - \sum\limits_{k=1}^{r} \lambda_{i_k} \omega'_{i_k} \wedge \omega_S {+}$$$$ \sum\limits_{j \notin S} \lambda_j
\sum\limits_{k=1}^{r} (-1)^{r+k} \omega'_{i_k} \wedge \omega_{(S_k,j)}=\sum\limits_{k=1}^{r} (-1)^{k+1} (d_x d_t + d_t d_x) \ln \alpha_{i_k} \wedge \omega_{S_k} + \left(\sum\limits_{j \notin S}
\lambda_j \omega'_j\right) \wedge \omega_T {+}$$$$ \sum\limits_{j \notin S} \lambda_j \sum\limits_{k=1}^{r} (-1)^{r+k} \omega'_{i_k} \wedge \omega_{(S_k,j)}=\left(\sum\limits_{j \notin S}
\lambda_j \omega'_j\right) \wedge \omega_S +\sum\limits_{j \notin S} \lambda_j \sum\limits_{k=1}^{r} (-1)^{r+k} \omega'_{i_k} \wedge \omega_{(S_k,j)}=(-1)^{r} \sum\limits_{j \notin S} \lambda_j
\omega'_{(S,j)}.$$

\noindent Hence, in the highest cohomology $H^{m}(\Omega(*\mathcal{A}),\nabla)$ we have

$$\nabla' \omega_S {\sim} (-1)^{l+1} \sum\limits_{j \notin S} \lambda_j \sum\limits_{k=1}^{l+1}(-1)^{k} d_x \ln \frac{\Delta_{(S,j)}}{\Delta_{((S,j)_k,n+1)}} \wedge \omega_{(S,j)_k}.
\eqno(4.3)$$

\noindent Now{, using (4.3),} we can represent $\nabla' \theta(S)$ as a linear combination (with {rational differential forms on $\mathbb{C}^n$ as coefficients}) of the sections $\omega_S$ of the vector bundle $\mathcal{H}$. Further it is necessary to expand the sections $\omega_S$ in the $\beta \text{nbc}$-basis, but {a} general formula for this purpose is {un}known. Below the corresponding computations for the integrals $J_2, I_n$ are given.\\

\begin{nscenter}\textbf{\hspace{0.5mm}5. Fuchsian system {for} $J_2(a, b, c, g; x, y)$}\end{nscenter}

In this section we introduce the notation $x=x_1, y=x_2, t=t_1, s=t_2$. The integrand in the integral $J_2(a, b, c, g; x, y)$ at the point $(x_0, y_0)$ such that $x_0 \neq 0, 1; y_0 \neq 0, 1; x_0 \neq y_0$ defines the following hyperplane
arrangement{:}
$$\mathcal{J}_2=\Big\{\{t=0\},\{t-1=0\},\{t-x_0=0\},\{s=0\},\{s-1=0\},\{s-y_0=0\},\{t-s=0\}\Big\}. \eqno(5.1)$$

\begin{center}\includegraphics[scale=0.44]{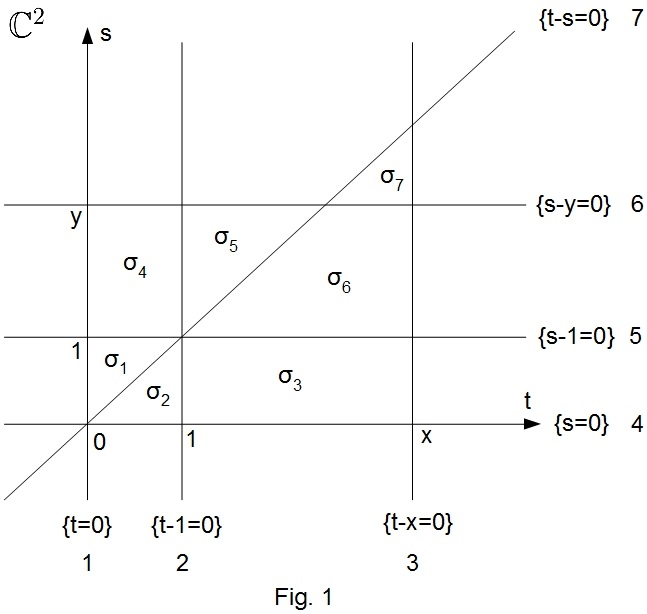}\end{center}

\noindent {We number the lines in this arrangement from 1 to 7 as shown on the figure.} Then $\beta nbc(\mathcal{J}_2)=\Big\{\{3, 5\}, \{2, 6\}, \{3, 7\}, \{6, 7\}, \{3, 6\}, {\{2, 7\}}, \{2, 5\}\Big\}$ is {a} $\beta\text{nbc}$-subset {for} the simplicial complex $\text{nbc}(\mathcal{J}_2)$. According to the formula (3.3) the corresponding collection of differential forms on $X=\mathbb{C}^2 \setminus \Big(\{x=0\} \cup \{x=1\} \cup \{y=0\} \cup \{y=1\} \cup \{x=y\}\Big)$, defining {a} $\beta\text{nbc}$-basis of the vector bundle $\mathcal{H}$, has the form

$$\eta_1= bc\frac{dt\wedge ds}{(t-x)(s-1)},~~~~ \eta_2=bc\frac{dt\wedge ds}{(t-1)(s-y)},$$

$$\eta_3=-cg\frac{dt\wedge ds}{(t-x)(t-s)},~~~~ \eta_4=-cg\frac{dt\wedge ds}{(s-y)(t-s)},~~~~\eta_5=c^2 \frac{dt\wedge ds}{(t-x)(s-y)}, \eqno(5.2)$$

$$\eta_6=-bg\frac{dt\wedge ds}{(t-1)(t-s)} - bg \frac{dt\wedge ds}{(s-1)(t-s)},~~~~\eta_7=b^2 \frac{dt\wedge ds}{(t-1)(s-1)} + bg \frac{dt\wedge ds}{(s-1)(t-s)}.$$

\noindent Using the identities $\nabla_{(x,y)}(\alpha){\sim} 0$ for the differential 1-forms $\alpha=d \ln{(t-x)},~ d \ln{(s-y)},~d \ln{(t-1)},~d \ln{(s-1)},\\d \ln{t} - d \ln{s},~ d \ln{(t-1)} - d \ln{(s-1)},~d \ln{t},~d \ln{s}$ one can obtain the following expansion of the sections $\omega_{ij}$ in the $\beta\text{nbc}$-basis

$$\omega_{34}{\sim}-\frac{1}{ac}(\eta_1+\eta_3+\eta_5),~~ \omega_{16}{\sim}-\frac{1}{ac}(\eta_2-\eta_4+\eta_5),~~ \omega_{24}{\sim}-\frac{1}{ab}(\eta_2+\eta_6+\eta_7),$$

$$\omega_{15}{\sim}-\frac{1}{ab}(\eta_1+\eta_7),~~ \omega_{25}{\sim}\frac{1}{b(2b+g)}(\eta_6+2\eta_7),~~ \omega_{27}{\sim}\frac{1}{bg}(\eta_6+\eta_7)-\frac{1}{g(2b+g)}(\eta_6+2\eta_7),$$

$$\omega_{14}{\sim}\frac{1}{a(2a+g)}(2\eta_1+2\eta_2+\eta_3-\eta_4+2\eta_5+\eta_6+2\eta_7),~~ \omega_{57}{\sim}-\frac{1}{bg} \eta_7 + \frac{1}{g(2b+g)}(\eta_6+2\eta_7),$$

$$\omega_{17}{\sim}-\frac{1}{g(2a+g)}(2\eta_1+2\eta_2+\eta_3-\eta_4+2\eta_5+\eta_6+2\eta_7)+\frac{1}{ag}(\eta_1+\eta_2-\eta_4+\eta_5+\eta_7),$$

$$\omega_{47}{\sim}\frac{1}{g(2a+g)}(2\eta_1+2\eta_2+\eta_3-\eta_4+2\eta_5+\eta_6+2\eta_7)-\frac{1}{ag}(\eta_1+\eta_2+\eta_3+\eta_5+\eta_6+\eta_7). $$

\noindent The corresponding matrix $T$, which is used in the formula (4.3), has the following form

$$T=\begin{pmatrix}
 0 & -1 & -x  & 0 & -1 & -y &  0 &  1  \\
 1 & 1 &   1 &  0 & 0 &   0 &  1 &   0 \\
 0 & 0 &   0 &  1 & 1 &   1 & -1 &  0 \\
\end{pmatrix}.$$

\noindent Further using the formula (4.3) it is easy to compute the expansions of $\nabla' \eta_i$ in the $\beta\text{nbc}$-basis (5.2), and after taking the integrals we obtain to the following proposition.

\vspace{1mm}\noindent\textbf{Proposition 1.} The collection of the vector-valued functions $f_i=\left(\int\limits_{\sigma_i}\Phi\eta_1, ..., \int\limits_{\sigma_i}\Phi\eta_7\right)^{T}, i=1,..,7$ is the basis of solutions of the Fuchsian system

$$d'f=\left(A\frac{d'x}{x}+B\frac{d'y}{y}+C\frac{d'(x-1)}{x-1}+D\frac{d'(y-1)}{y-1}+E\frac{d'(x-y)}{x-y}\right)f,~ \text{where}$$

$$A=\begin{pmatrix}
 a+c & 0 & 0 & 0 & 0 & 0 & c  \\
 0 & 0 & 0 & 0 & 0 & 0 & 0 \\
 g & 0 & 2a+c+g & c & g & c & 0 \\
 0 & 0 & 0 & 0 & 0 & 0 & 0 \\
 0 & c & 0 & -c & a+c & 0 & 0 \\
 0 & 0 & 0 & 0 & 0 & 0 & 0 \\
 0 & 0 & 0 & 0 & 0 & 0 & 0
\end{pmatrix},~~~~
B=\begin{pmatrix}
 0 & 0 & 0 & 0 & 0 & 0 & 0  \\
 0 & a+c & 0 & 0 & 0 & c & c \\
 0 & 0 & 0 & 0 & 0 & 0 & 0 \\
 0 & -g & c & 2a+c+g & -g & c & 0 \\
 c & 0 & c & 0 & a+c & 0 & 0 \\
 0 & 0 & 0 & 0 & 0 & 0 & 0 \\
 0 & 0 & 0 & 0 & 0 & 0 & 0
\end{pmatrix},$$

$$C=\begin{pmatrix}
 b+g & 0 & -b & 0 & 0 & 0 & -c  \\
 0 & c & 0 & 0 & -b & 0 & 0 \\
 -g & 0 & 2b & 0 & 0 & -c & 0 \\
 0 & 0 & 0 & 0 & 0 & 0 & 0 \\
 0 & -c & 0 & 0 & b & 0 & 0 \\
 g & 0 & -2b & 0 & 0 & c & 0 \\
 -b-g & 0 & b & 0 & 0 & 0 & c
\end{pmatrix},~~~~
D=\begin{pmatrix}
 c & 0 & 0 & 0 & -b & 0 & 0  \\
 0 & b+g & 0 & b & 0 & -c & -c \\
 0 & 0 & 0 & 0 & 0 & 0 & 0 \\
 0 & g & 0 & 2b & 0 & -c & 0 \\
 -c & 0 & 0 & 0 & b & 0 & 0 \\
 0 & -g & 0 & -2b & 0 & c & 0 \\
 0 & -b & 0 & b & 0 & 0 & c
\end{pmatrix},$$

$$E=\begin{pmatrix}
 0 & 0 & 0 & 0 & 0 & 0 & 0  \\
 0 & 0 & 0 & 0 & 0 & 0 & 0 \\
 0 & 0 & c & -c & -g & 0 & 0 \\
 0 & 0 & -c & c & g & 0 & 0 \\
 0 & 0 & -c & c & g & 0 & 0 \\
 0 & 0 & 0 & 0 & 0 & 0 & 0 \\
 0 & 0 & 0 & 0 & 0 & 0 & 0
\end{pmatrix}.$$

\begin{nscenter}\textbf{\hspace{0.5mm}6. Fuchsian system {for} $I_n(a_0, a_1, ..., a_n, g; x_2,...,x_n)$}\end{nscenter}

In this section we introduce the notation $x_0=x_0^0=0, x_1=x_1^0=1$. The integrand in the integral $I_n(a_0,a_1,...,a_n,g; x_2,...x_n)$ at the point $x_i^0 \neq x_j^0,~ 0 \leq i \neq j \leq n$ defines the following hyperplane arrangement

$$\mathcal{I}_n=\Big\{V_i,H_i,D ~ | ~ i=0,..,n \Big\},~~ V_i=\{t-x^{0}_i=0\},~ H_i=\{s-x^{0}_i=0\}, ~ D=\{t-s=0\}. \eqno(6.1)$$

\begin{center}\includegraphics[scale=0.44]{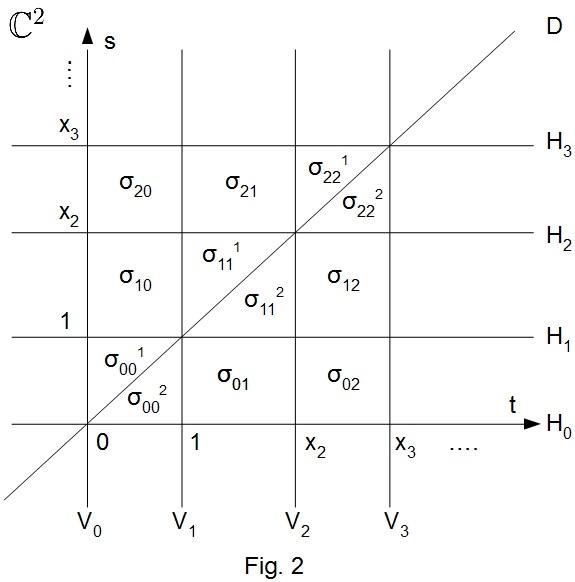}\end{center}

\noindent It easy to see that the $\beta \text{nbc}$-subset of the simplicial complex $\text{nbc}(\mathcal{I}_n)$ is $\beta \text{nbc}(\mathcal{I}_n)=\Big\{\{V_i,H_j\} ~ | ~ 1 \leq i \leq n, 1 \leq j \leq n \Big\} \cup \Big\{\{V_i,D\} ~ | ~ 1 \leq i \leq n \Big\}$. Now using the formula (3.3) we compute the $\beta\text{nbc}$-basis of the vector bundle $\mathcal{H}$ which is defined by the following differential forms
$$\theta(V_i,H_j)=a_i a_j V_i H_j, ~~~ \theta(V_i, H_i)= a_i^2 V_i H_i+ a_i g D H_i, ~~~ \theta(V_i, D)=a_i g V_i D + a_i g H_i D, ~~ 1 \leq i \neq j \leq n. \eqno(6.2)$$

\noindent Here we use the notation $MN:=d\log \alpha_M \wedge d\log \alpha_N$ for the hyperplanes $M,N$ from {the} arrangement $\mathcal{I}_n$, where $\alpha_M,\alpha_N$ are the linear polynomials defining the hyperplanes $M,N$.

The following identities are necessary {in order} to obtain the expansions of the sections $MN$ corresponding {to} all pairs of hyperplanes from $\mathcal{I}_n$ in the $\beta\text{nbc}$-basis (6.2){:}

$$\nabla_x(H_k)=\sum\limits_{i=0}^{n} a_i V_i H_k + g D H_k = a_0 V_0 H_k + \sum\limits_{i=1}^{n} \frac{1}{a_k} \theta(V_i, H_k) {\sim} 0, ~~ 1 \leq k \leq n;$$

$$\nabla_x(V_k)=-\sum\limits_{i=1}^{n} a_i V_k H_i - g V_k D = -a_0 V_k H_0 - \frac{1}{a_k}\Big(\theta(V_k,H_k)+\theta(V_k,D)\Big)-\sum\limits_{i \neq k} \frac{1}{a_k} \theta(V_k, H_i) {\sim}
0;$$

$$\nabla_x(H_k-V_k)=\sum a_i V_i H_k - \sum a_i H_i V_k + g D(H_k-V_k)=(2a_k+g)V_k H_k + \sum\limits_{i \neq k}a_i V_i H_k - \sum\limits_{i \neq k} a_i H_i V_k {=}$$$$ (2a_k+g)V_k H_k + a_0(V_0
H_k - H_0 V_k) + \sum\limits_{i \neq 0,k}a_i(V_i H_k - H_i V_k){=}$$$$(2a_k+g)V_k H_k - \frac{1}{a_k}\left(\theta(V_k,D)+\sum\limits_{i \neq 0}\theta(V_i,H_k)+\sum\limits_{i \neq
0}\theta(V_k,H_i)\right)+\sum\limits_{i \neq 0,k} \frac{1}{a_k} \Big(\theta(V_i,H_k)+\theta(V_k,H_i)\Big) {\sim} 0;$$

$$\nabla_x(H_0-V_0)=\sum a_i V_i H_0 - \sum a_i H_i V_0 + g D(H_0-V_0)=(2a_0+g) V_0 H_0 + \sum\limits_{i \neq 0} a_i (V_i H_0 - H_i V_0){=}$$$$(2a_0+g)V_0 H_0 - \sum\limits_{i \neq
0}\frac{1}{a_0}\left(\theta(V_i,D)+\sum\limits_{j \neq 0} \theta(V_i,H_j)+\sum\limits_{j \neq 0}\theta(V_j,H_i)\right) {\sim} 0;$$

$$\nabla_x(V_0)=a_0 H_0 V_0 + \sum\limits_{i \neq 0} a_i H_i V_0 + g D V_0 {=}$$$$ g D V_0 - \frac{1}{2a_0+g}\left(\sum\limits_{i \neq 0}\theta(V_i,D)+2\sum\limits_{i,j \neq
0}\theta(V_i,H_j)\right)+\frac{1}{a_0}\sum\limits_{i,j \neq 0} \theta(V_j,H_i) {\sim} 0;$$

$$\nabla_x(H_0)=\sum a_i V_i H_0 + g D H_0 {=}$$$$ g D H_0 + \frac{1}{2a_0+g}\left(\sum\limits_{i \neq 0}\theta(V_i,D)+2\sum\limits_{i,j \neq 0}\theta(V_i,H_j)\right)-\frac{1}{a_0}\sum_{i \neq 0}
\Big(\theta(V_i,D)+\sum\limits_{j \neq 0}\theta(V_i,H_j)\Big) {\sim} 0;$$

$$\theta(V_k,H_k)=a_k^2 V_k H_k + a_k g D H_k = a_k g D H_k + \frac{a_k}{2a_k+g}\Big(\theta(V_k,D)+2\theta(V_k,H_k)\Big);$$

$$\theta(V_k,D)=a_k g V_k D + a_k g H_k D = a_k g V_k D + \frac{a_k}{2a_k+g}\Big(\theta(V_k,D)+2\theta(V_k,H_k)\Big)-\theta(V_k,H_k).$$

\noindent {The following expansions are straightforward consequences of the above identities:}

$$V_0 H_k {\sim} - \frac{1}{a_0 a_k} \sum\limits_{i \neq 0} \theta(V_i,H_k),~~ V_k H_0 {\sim} -\frac{1}{a_0 a_k}\left(\theta(V_k,D)+\sum\limits_{i \neq 0}\theta(V_k,H_i)\right),$$

$$V_0 H_0 {\sim} \frac{1}{(2 a_0 + g)a_0}\left(\sum\limits_{i \neq 0}\theta(V_i,D)+2\sum\limits_{i,j \neq 0}\theta(V_i,H_j)\right),~~ V_k H_k {\sim} \frac{1}{(2 a_k + g) a_k} \Bigg( \theta(V_k,
D)+2\theta(V_k,H_k) \Bigg),$$

$$V_0 D {\sim} - \frac{1}{(2a_0+g)g}\left(\sum\limits_{i \neq 0}\theta(V_i,D)+2\sum\limits_{i,j \neq 0}\theta(V_i, H_j)\right)+\frac{1}{a_0 g} \sum\limits_{i,j \neq 0}\theta(V_i,H_j),
\eqno(6.3)$$

$$V_k D {\sim} - \frac{1}{(2a_k+g)g}\Bigg(\theta(V_k,D)+2\theta(V_k,H_k)\Bigg)+\frac{1}{a_k g}\Bigg(\theta(V_k,D)+\theta(V_k,H_k)\Bigg),$$

$$H_0 D {\sim} \frac{1}{(2 a_0 +g)g}\left(\sum\limits_{i \neq 0} \theta(V_i,D)+2 \sum\limits_{i, j \neq 0}\theta(V_i,H_j)\right)-\frac{1}{a_0 g}\left(\sum\limits_{i \neq
0}\theta(V_i,D)+\sum\limits_{i,j \neq 0}\theta(V_i,H_j)\right);$$

$$ H_k D {\sim} \frac{1}{(2a_k+g)g}\Bigg(\theta(V_k, D)+2 \theta(V_k,H_k)\Bigg) - \frac{1}{a_k g}\theta(V_k,H_k), ~~ \text{где} ~~ 1 \leq k \leq n.$$

\noindent {It is useful} to note that all {the} determinants appearing in the formula (4.5) have the following simple form{:} $\Delta_{(V_i,H_j,V_k)}=x_i-x_k,~ \Delta_{(V_i,H_j,H_k)}=x_j-x_k,~ \Delta_{(V_i,H_j,D)}=x_i-x_j,~ \Delta_{(V_i,V_j,D)}=x_i-x_j,~ \Delta_{(H_i,H_j,D)}=x_i-x_j$. Further using the formula (4.5){,} one can {compute:}
$$\nabla'(V_i H_j) {\sim} \frac{d'(x_i-x_j)}{x_i-x_j}\wedge \Bigg(a_j (H_j V_j+V_i H_j) + a_i (H_i V_i +V_i H_j) + g (H_j D - V_i D + V_i H_j) \Bigg) {+} $$$$ \sum\limits_{k \neq j}
\frac{d(x_i-x_k)}{x_i-x_k} \wedge a_k (H_j V_k + V_i H_j) + \sum\limits_{k \neq i} \frac{d'(x_j-x_k)}{x_j-x_k} \wedge a_k (H_k V_i + V_i H_j),$$

$$\nabla'(V_i D) {\sim} \sum\limits_{k=0}^n \frac{d'(x_i-x_k)}{x_i-x_k} \wedge a_k (D V_k + D H_k + 2 V_i D + H_k V_i), \eqno(6.4)$$

$$\nabla'(H_i D) {\sim} \sum\limits_{k=0}^n \frac{d'(x_i-x_k)}{x_i-x_k} \wedge a_k (D V_k + D H_k + 2 H_i D + V_k H_i).$$

\noindent Substituting (6.4) into (6.6) we obtain the expansions of $\nabla'(V_i H_j),~ \nabla'(V_i D),~ \nabla'(H_i D)$ in the $\beta\text{nbc}$-basis (6.2). Now we immediately obtain the connection matrix in the $\beta\text{nbc}$-basis. Therefore, after taking the integrals the following proposition {follows}.

\vspace{1mm}\noindent\textbf{Proposition 2.} The collection of the vector-valued functions $f^{rsp}=\Big(f^{rsp}_{ij},f^{rsp}_k\Big)_{i,j,k=1}^{n}$, where $f^{rsp}_{ij}=\int\limits_{\sigma_{rs}^{p}} \theta(V_i,H_j),~ f^{rsp}_k=\int\limits_{\sigma_{rs}^{p}}
\theta(V_k,D),~ 0 \leq r, s \leq n-1 ~\text{and}~ p=1,2 ~\text{for}~ r=s${,} is the basis of solutions of the following Fuchsian system (the indices $r,s,p$ are omitted){:}

$$d'f_{ij}=\frac{d'(x_i-x_j)}{x_i-x_j}\wedge\bigg(\big(a_i+a_j+g\big)f_{ij}-a_i f_{jj} - a_j f_{ii} - a_j f_{i}\bigg) + \frac{d'x_i}{x_i} \wedge \bigg(a_i \sum\limits_{l \neq 0} f_{lj} + a_0
f_{ij}\bigg) +$$$$+ \sum\limits_{k \neq 0,j} \frac{d'(x_i-x_k)}{x_i-x_k} \wedge \bigg(a_k f_{ij}-a_i f_{kj}\bigg) + \frac{d'x_j}{x_j} \wedge \bigg(a_j \Big(f_{i}+\sum\limits_{l \neq 0}
f_{il}\Big)+ a_0 f_{ij}\bigg) + \sum\limits_{k \neq 0,i} \frac{d'(x_j-x_k)}{x_j-x_k} \wedge \bigg( a_k f_{ij}-a_j f_{ik}\bigg),$$

$$d'f_{ii}=\frac{d'x_i}{x_i}\wedge\bigg( \sum\limits_{l \neq 0} \Big( (a_i+g ) f_{li} + a_i ( f_{il} - f_{l}) \Big) + 2a_0 f_{ii} + a_i f_i \bigg)+\sum\limits_{k \neq 0,i}
\frac{d'(x_i-x_k)}{x_i-x_k} \wedge \bigg(a_i \Big(f_{k}-f_{ik}-f_{ki}\Big)+2 a_k f_{ii}-g f_{ki}\bigg),$$

$$d'f_{i}=\frac{d'x_i}{x_i}\wedge\bigg( (2a_0+g) f_{i}+ \sum\limits_{l \neq 0} \Big( 2 a_i f_{l}+g ( f_{il} - f_{li} ) \Big) \bigg)+\sum\limits_{k \neq 0, i}\frac{d'(x_i-x_k)}{x_i-x_k} \wedge
\bigg(2 a_k f_{i}-2 a_i f_{k}+g\Big(f_{ki}-f_{ik}\Big)\bigg).$$

\vspace{3mm}\begin{nscenter}\textbf{\hspace{0.5mm}7. {S}ystem of partial differential equations for $J_n(a,b,c,g;x_1,...,x_n)$}\end{nscenter}

\vspace{1mm} It is easy to see that the sections of the vector bundle $\mathcal{H}$ defined by the differential forms $\prod\limits_{i=1}^{n}(t_i-x_i)^{-k_i}dt$ correspond to the partial derivatives of the integral $J_n(a,b,c,g;x_1,...,x_n)$. {To obtain a} system of differential equations for the integral $J_n(a,b,c,g;x_1,...,x_n)${, one has to use the relations $\nabla_x (\beta_k)=0$,} where

$$\beta_k=\sum\limits_{i=1}^{n}(-1)^{i+1} \frac{t_i(t_i-1)}{t_i-x_i}\prod\limits_{j\neq i}\left(1+\frac{x_j-x_k}{t_j-x_j}\right)dt_1\wedge ... \wedge dt_{i-1} \wedge dt_{i+1} \wedge ... \wedge
dt_n , ~~ k=1,..,n. \eqno(7.1)$$

\noindent The explicit form of these relations is given below{:}

$$0=\sum\limits_{i=1}^{n}\left(\frac{2t_i-1}{t_i-x_i}-\frac{t_i(t_i-1)}{(t_i-x_i)^2}\right)\prod\limits_{j\neq i}\left(1+\frac{x_j-x_k}{t_j-x_j}\right)dt{+}$$
$$\sum\limits_{i=1}^{n}\left(a\frac{t_i-1}{t_i-x_i}+b\frac{t_i}{t_i-x_i}+c\frac{t_i(t_i-1)}{(t_i-x_i)^2}\right)\prod\limits_{j\neq i}\left(1+\frac{x_j-x_k}{t_j-x_j}\right)dt{+}$$
$$g\sum\limits_{i<j}\frac{1}{\prod\limits_{s=1}^{n}(t_s-x_s)}\frac{1}{t_i-t_j}\left(t_i(t_i-1)\prod\limits_{s\neq i}(t_s-x_k)-t_j(t_j-1)\prod\limits_{s\neq j}(t_s-x_k)\right)dt{=}$$
$$\sum\limits_{i=1}^{n}\prod\limits_{j\neq i}\left(1+\frac{x_j-x_k}{t_j-x_j}\right) \cdot
\left(a\frac{t_i-1}{t_i-x_i}+\frac{2t_i-1}{t_i-x_i}+b\frac{t_i}{t_i-x_i}+(c-1)\frac{t_i(t_i-1)}{(t_i-x_i)^2}\right)dt{+}$$
$$g\sum\limits_{i<j}\frac{\prod\limits_{s\neq i,j}(t_s-x_k)}{\prod\limits_{s=1}^{n}(t_s-x_s)} \cdot \frac{1}{t_i-t_j} \cdot \left( t_i(t_i-1)(t_j-x_k) - t_j(t_j-1)(t_i-x_k) \right) dt.
\eqno(7.2)$$

\noindent Now we {apply} some transformations{:}

$$a\frac{t_i-1}{t_i-x_i}=a+a\frac{x_i-1}{t_i-x_i}, ~~ \frac{2t_i-1}{t_i-x_i}=2+\frac{2x_i-1}{t_i-x_i}, ~~ b\frac{t_i}{t_i-x_i}=b+b\frac{x_i}{t_i-x_i},$$
$$(c-1)\frac{t_i(t_i-1)}{(t_i-x_i)^2}=(c-1)\left(1+\frac{2x_i-1}{t_i-x_i}+\frac{x_i(x_i-1)}{(t_i-x_i)^2}\right),$$
$$t_i(t_i-1)(t_j-x_k)-t_j(t_j-1)(t_j-x_k)=t_{i}^{2}t_j-t_{i}^{2}x_k-t_i t_j - t_{j}^{2}t_i+t_{j}^{2}x_k+t_j t_i - t_j x_k=$$
$$=t_i t_j (t_i-t_j)-(t_i-t_j)(t_i+t_j)x_k+(t_i-t_j)x_k=(t_i-t_j)(t_i t_j - (t+i+t_j-1)x_k),$$
$$\frac{\prod\limits_{s\neq i,j}(t_s-x_k)}{\prod\limits_{s=1}^{n}(t_s-x_s)}=\prod\limits_{s\neq i, j}\left( 1+\frac{x_s-x_k}{t_s-x_s}\right) \cdot \frac{1}{(t_i-x_i)(t_j-x_j)}.$$

\noindent {Plugging them in} the formulas (7.2){, we obtain:}
$$0=\sum\limits_{i=1}^{n}\prod\limits_{j\neq i}\left(1+\frac{x_j-x_k}{t_j-x_j}\right) \cdot \left( (1+a+b+c)+\frac{(a+b+2c)x_i-(a+c)}{t_i-x_i} + (c-1) \frac{x_i(x_i-1)}{(t_i-x_i)^2}
\right)dt{+}$$
$$g\sum\limits_{i<j}\prod\limits_{s\neq i, j}\left( 1+\frac{x_s-x_k}{t_s-x_s}\right) \cdot \frac{t_i t_j - (t_i+t_j-1)x_k}{(t_i-x_i)(t_j-x_j)}dt. \eqno(7.3)$$

\noindent The last fraction of the identity (7.3) {expands} over the differential forms $\prod\limits_{i=1}^{n}(t_i-x_i)^{-k_i}dt$ corresponding {to} the partial derivatives of {our} integral{:}
$$\frac{t_i t_j - (t_i+t_j-1)x_k}{(t_i-x_i)(t_j-x_j)}dt=1+\frac{x_i}{t_i-x_i}+\frac{x_j}{t_j-x_j}+\frac{x_i
x_j}{(t_i-x_i)(t_j-x_j)}-\frac{x_k}{t_i-x_i}-\frac{x_k}{t_j-x_j}-\frac{(x_i+x_j-1)x_k}{(t_i-x_i)(t_j-x_j)}=$$
$$=1+\frac{x_i-x_k}{t_i-x_i}+\frac{t_j-x_k}{t_j-x_j}+\frac{x_i x_j-(x_i+x_j-1)x_k}{(t_i-x_i)(t_j-x_j)}.$$

\noindent Now after replacing of the differential forms by the corresponding partial derivatives we obtain the following proposition.

\vspace{1mm}\noindent\textbf{Proposition 3.} The integral $J_n(a,b,c,g;x_1,...,x_n)$ is the solution of the system of $n$ partial differential equations of the { order $(n+1)$:}

$$P_k J_n=0, ~~ 1 \le k \le n, ~\text{where}$$
$$P_k=\sum\limits_{i=1}^{n}\prod\limits_{j\neq i}\left(c-(x_j-x_k)\partial_j\right) \cdot \left( (1+a+b+c)c-((a+b+2c)x_i-(a+c))\partial_i + x_i(x_i-1)\partial_i^{2} \right){+}$$
$$g\sum\limits_{i<j}\prod\limits_{s\neq i, j}\left( c-(x_s-x_k)\partial_s\right) \cdot \left( c^2-c(x_i-x_k)\partial_i-c(x_j-x_k)\partial_j+(x_i
x_j-(x_i+x_j-1)x_k)\partial_i\partial_j\right).$$

In particular, for $n=2$ this system has the form ($x=x_1, y=x_2, u(x,y)=J_2(a,b,c,g;x,y)$){:}
\begin{equation*}
 \begin{cases}
 x(1-x)(x-y)u_{xxy}+cx(1-x)u_{xx}+cy(1-y)u_{yy}{+}\\
 ((x-y)((a+b+2c)x-(a+c))-gx(1-x))u_{xy}{+}\\
 ((a+b+2c)x-(a+c))cu_x-((x-y)(1+a+b+c+g){-}\\
 (a+b+2c)y+(a+c))cu_y-(2(1+a+b+c)+g)c^{2}u=0,\\
 \\
 y(1-y)(y-x)u_{xyy}+cx(1-x)u_{xx}+cy(1-y)u_{yy}{+}\\
 ((y-x)((a+b+2c)y-(a+c))-gy(1-y))u_{xy}{+}\\
 ((a+b+2c)y-(a+c))cu_y-((y-x)(1+a+b+c+g){-}\\
 (a+b+2c)x+(a+c))cu_x-(2(1+a+b+c)+g)c^{2}u=0.
 \end{cases}
\end{equation*}

\begin{nscenter}\textbf{\hspace{0.5mm}8. {Another} Fuchsian system {for} $I_2(a,b,c,g;z)$}\end{nscenter}

\vspace{1mm} It is known that the Dotsenko{--}Fateev integral $I_2(a, b, c, g; z)$ {satisfies} the following ordinary differential equation of third order [1, 2]{:}
$$y'''+\frac{K_1z+K_2(z-1)}{z(z-1)}y''+\frac{L_1z^2+L_2(z-1)^2+L_3z(z-1)}{z^2(z-1)^2}y'+\frac{M_1z+M_2(z-1)}{z^2(z-1)^2}y=0, \eqno(8.1)$$
\noindent where coefficients have the form

$$K_1=-(3b+3c+g),~~ K_2=-(3a+3c+g),$$
$$L_1=(b+c)(1+2b+2c+g),~~ L_2=(a+c)(1+2a+2c+g),~~  \eqno(8.2)$$
$$L_3=(b+c)(1+2a+2c+g)+(a+c)(1+2b+2c+g)+(c-1)(a+b+c)+(3c+g)(1+a+b+c+g),$$
$$M_1=-c(2+2a+2b+2c+g)(1+2b+2c+g),~~ M_2=-c(2+2a+2b+2c+g)(1+2a+2c+g).$$

Denote $f_0=(y,y',y'')^T$ and write (8.1) as {the} Pfaffian system $\frac{df_0}{dz}=\Omega_0 f_0$, where
$$
\Omega_0=\begin{pmatrix}
 0 & 1 & 0 \\
 0 & 0 & 1 \\
 -\frac{M_1}{z(z-1)^2}-\frac{M_2}{z^2(z-1)} & -\frac{L_1}{(z-1)^2}-\frac{L_2}{z^2}-\frac{L_3}{z(z-1)} & -\frac{K_1}{z-1}-\frac{K_2}{z}
\end{pmatrix}.\eqno(8.3)
$$

\noindent Now transform this system{, applying the} multiplication of the vector-valued function $f_0$ {by} the following matrix{:}
$$
\Gamma_0=\begin{pmatrix}
 1 & 0 & 0 \\
 0 & z-1 & 0 \\
 0 & 0 & (z-1)^2
\end{pmatrix},~~
f_1=\Gamma_0 f_0.  \eqno(8.4)
$$

\noindent Then the system (8.3) takes the form $\frac{d}{dz} f_1= \Omega_1 f_1$, where
$$\Omega_1=\Gamma_0 \Omega_0 \Gamma_0^{-1}-\Gamma_0 \frac{d}{dz} \Gamma_0^{-1}=
\begin{pmatrix}
 0 & \frac{1}{z-1} & 0 \\
 0 & \frac{1}{z-1} & \frac{1}{z-1} \\
 -\frac{M_1+M_2}{z}+\frac{M_2}{z^2} & -\frac{L_2+L_3}{z}-\frac{L_1}{z-1}+\frac{L_2}{z^2} & -\frac{K_2}{z}+\frac{2-K_1}{z-1}
\end{pmatrix}. \eqno(8.5)$$

\noindent The obstruction {to} being Fuchsian {for} the system (8.5) is the presence of summands of the form $\frac{C}{z^2}$ {in} the lower row of the matrix $\Omega_1$. {To remove} these summands{,} we transform the system (8.5){, applying the multiplication} of the vector-valued function $f_1$ {by} the following lower triangular matrix{:}

$$\Gamma_1=
\begin{pmatrix}
 1 & 0 & 0 \\
 0 & 1 & 0 \\
  \frac{\eta}{z} & \frac{\zeta}{z} & 1
\end{pmatrix},~~ f_2=\Gamma_1 f_1{.} \eqno(8.6)
$$

\noindent {The coefficients $\eta, \zeta$ are chosen to} satisfy {the} two equations $M_2 +\eta \zeta + \eta K_2 - \eta = 0, ~ L_2+\zeta K_2 + \zeta^2 - \zeta = 0${,} which provide {the} vanishing of {the} summands of the form $\frac{C}{z^2}$ in the lower row of the matrix $\Omega_2$ of the {transformed} system (it is easy to see that this system has two solutions for $L_2 \neq 0$ and one solution for $L_2=0$). {The resulting} matrix $\Omega_2$ has the following form{:}
$$
\Omega_2=\Gamma_1 \Omega_1 \Gamma_1^{-1} - \Gamma_1 \frac{d}{dz} \Gamma_1^{-1}=\eqno(8.8)
$$
$$=
\begin{pmatrix}
 0 & \frac{1}{z-1} & 0 \\
 -\frac{\eta}{z(z-1)} & \frac{1}{z-1}-\frac{\zeta}{z(z-1)} & \frac{1}{z-1} \\
 -\frac{M_1+M_2}{z}-\frac{\eta \zeta}{z-1}+\frac{\eta \zeta}{z} -\frac{(2-K_1) \eta}{z(z-1)} & \frac{\eta + \zeta}{z(z-1)}-\frac{L_2+L_3}{z}-\frac{L_1}{z-1}-\frac{\zeta^2}{z-1}+\frac{\zeta^2}{z}
 - \frac{(2-K_1)\zeta}{z(z-1)} & \frac{\zeta}{z(z-1)}-\frac{K_2}{z}+\frac{2-K_1}{z-1}
 \end{pmatrix}
$$

\noindent {By collecting similar terms,} we obtain the following proposition.

\noindent\textbf{Proposition 4.} Let $y$ be {a} solution of the differential equation (8.1){, and let the} complex numbers $\eta, \zeta$ satisfy {the} two equations $M_2 +\eta \zeta + \eta K_2 - \eta = 0, ~ L_2+\zeta K_2 + \zeta^2 - \zeta = 0$. Then the vector-valued function $f=\Gamma_1 \Gamma_0 f_0=(y,(z-1)
y',\frac{\eta}{z}y+\zeta\frac{z-1}{z}y'+(z-1)^2 y'')^T$ is {a} solution of the Fuchsian system
$$df=\left(A\frac{dz}{z}+B\frac{dz}{z-1}\right)f,~ \text{where}$$
$$
A=
\begin{pmatrix}
 0 & 0 & 0 \\
 \eta & \zeta & 0 \\
  -M_1-M_2+(\zeta + 2-K_1)\eta & -\eta-L_2-L_3+\zeta^2+(1-K_1)\zeta & -\zeta-K_2
\end{pmatrix},
$$
$$
B=
\begin{pmatrix}
 0 & 1 & 0 \\
 -\eta & 1-\zeta & 1 \\
  (K_1-\zeta-2)\eta & \eta-L_1-\zeta^2-(1-K_1)\zeta & 2+\zeta-K_1
\end{pmatrix}.
$$
\noindent Moreover, the map $y \longmapsto f$ is {an} isomorphism between the space of solutions of (8.1) and the space of solutions of {the} Fuchsian system.

\vspace{0.5mm}\begin{nscenter}\textbf{Acknowledgements}\end{nscenter}

We {are pleased to acknowledge} useful {discussions with S.~Tanabe, Vl.~Dotsenko, Vl.~Fateev and A.~B.~Belavin}.

\vspace{5mm}\begin{nscenter}\textbf{References}\end{nscenter}

\vspace{1mm}\noindent 1. Vl. S. Dotsenko, Vl. A. Fateev, Conformal algebra and multipoint correlation functions in 2D statistical models. Nuclear Physics B, V. 240, №3, 1984.

\noindent 2. V. Toledano Laredo, Fusion of positive energy representations of LSpin(2n). Arxiv preprint math/0409044, 2004.

\noindent 3. P. Orlik, H. Terao, Arrangements and hypergeometric integrals.  MSJ Memoir, 9, Math. Soc. Japan, Tokyo, 2001.

\noindent 4. P. Appell, M. J. {Kamp\'e de F\'eriet}, Fonctions {hyperg\'eom\'etriques} et {hypersph\'eriques}. Paris, Gauthier -- Villars, 1926.

\noindent 5. V. A. Golubeva, The hypergeometric functions of two variables of Appell and {Kamp\'e de F\'eriet}, Siberian journal of mathematics, v. XX, N. 5.

\noindent 6. T. Oaku, Computation of the characteristic variety and the singular locus of a system of differential equations with polynomial coefficients. Japan J. Indust. Appl. Math. 11 (1994),
485-497.

\noindent 7. H. Esnault, V. Schechtman, E. Viehweg, Cohomology of local systems of the complement of hyperplanes, Invent. math. 109 (1992), 557-561; Erratum. ibid. 112 (1993), 447.

\noindent 8. M. Falk, H. Terao, $\beta\text{nbc}$-bases for cohomology of local systems on hyperplane complements, Trans. AMS, 349 (1997), 189-202.

\noindent 9. P. Orlik, H. Terao, Arrangements of hyperplanes, Springer Verlag, 1992.

\noindent 10. M. Kato, Connection formulas for Appell's system $F_4$ and some applications, Funkcialaj Ekvacioj, 38 (1995), 243-266.

\end{document}